\documentclass[12pt]{article}
\usepackage[cp1251]{inputenc}
\usepackage[ukrainian, russian, english]{babel}

\usepackage{ amsfonts, amssymb}
\usepackage[mathscr]{eucal}

\usepackage{mathrsfs}

\begin{document}
\selectlanguage{ukrainian} \thispagestyle{empty}
 \pagestyle{myheadings}               %%%%%%%%%%%%%%%%%%%%% <--------------
%%%%%%%%%%%%%%%%%%%%%%%%%%%%%%%%%%%%%%%%%%%%%%%%%%%%%%%%%%%%%%%%%%%%%%%%%%%%%
%%%%%%%%%%%%%%%%%%%%%%%%%%%%%%%%%%%%%%%%%%%%%(\ref {3})%%%%%%%%%%%%%%%%%%%%%%
%%%%%%%%%%%%%%%%%%%%%%%%%%%%%%%%%%%%%%%%%%%%%%%%%%%%%%%%%%%%%%%%%%%%%%%%%%%%%

\vskip 1cm

{\bf J.~Prestin{$^1$}, V.\,V.~Savchuk{$^2$}, A.\,L.~Shidlich{$^3$}

\vskip 0.5cm

Direct and inverse  approximation  theorems  

of  $2\pi$-periodic functions 
by Taylor--Abel--Poisson means}

\vskip 0.5cm

{\it \small Institute of Mathematics, University of L{\"{u}}beck$^1$}

{\it \small Institute of Mathematics of the National Academy of Sciences of Ukraine$^{2,3}$}

\vskip 0.5cm

{\small {\it  Emails:}  prestin@math.uni-luebeck.de$^1$,  vicsavchuk@gmail.com$^2$, shidlich@gmail.com$^3$.}

\vskip 0.5cm

{\it We obtain direct and inverse  approximation  theorems of  $2\pi$-periodic functions by Taylor--Abel--Poisson operators in the integral metric.}

\vskip 0.5cm


\vskip 0.5cm


{\it UDC:} 517.5

\vskip 1.5cm

It is well-known that  any function $f \in L_p $, $f \not \equiv {\rm const}$, can be approximated by its Abel-Poisson means $f(\varrho,\cdot)$ with a precision not better than $1-\varrho$. It relates to the so-called saturation property of this approximation method. From this property, it follows that for any $f\in L_p$, the relation $\|f-f(\varrho,\cdot)\|_p=o(1-\varrho)$, $\varrho\to 1-$, holds only in the trivial case where  $f\equiv{\rm const}$. Therefore, any additional restrictions on the smoothness of functions don't give us the order of approximation better than  $1-\varrho$. In this connection, a natural question is
%%%the finding of
to find
a linear operator,  constructed similarly to the
%%%operator Poisson
Poisson operator, which takes into account the smoothness properties of functions and at the same time, for a given functional class,  is the best in a certain sense.
In [\ref{Savchuk_Zastavnyi}], for classes of convolutions, whose kernels were generated by some moment sequences, the authors proposed the general method of construction of similar operators that take into account properties of such kernels
%%%(
and hence, the smoothness of functions from corresponding  classes.
%%%).
One example of such operators are the operators $ A_{\varrho, r} $, which are the main subject of study in this paper.

The operators $A_{\varrho,r}$ were first studied in [\ref{Savchuk}], where in the terms of these operators, the author gave the structural characteristic of Hardy-Lipschitz classes $H^r_p\mathop{\rm Lip}\alpha$  of
%%%one variable
functions of one variable,
holomorphic
%%%in
on
the unit circle
%%%in
of
the complex plane. In  [\ref{Savchuk_Shidlich}],  in
%the
terms of approximation estimates of such operators in
%%%the
some
spaces $S^p$ of Sobolev type, %%% is it ok to add "of Sobolev type" ??? It would be better readable.
the authors give a constructive description of classes of functions  of several variables, whose generalized derivatives belong to the classes $S^pH_\omega$. Similar operators of polynomial type were studied in [\ref{Chui}], [\ref{Holland}], [\ref{Mohapatra}], [\ref{Chandra}] etc. In particular, in [\ref{Chui}], the authors found the  degree of convergence of the well-known Euler and Taylor means to the functions $f$ from some subclasses of the Lipschitz classes ${\rm Lip} \alpha$ in the uniform norm. In [\ref{Mohapatra}], the analogical results for Taylor means were obtained in the $L_p$--norm.

In this paper, we continue the study of
%%%approximative
approximation
properties of the operators $A_{\varrho,r}$. In particular, we find the relation of these operators with the operators $L_{\varrho, r}$ and $B_{\varrho,r}$, considered in [\ref{Leis}] and [\ref{Butzer_Sunouchi}]. Also we give direct and inverse approximation theorems by the operators $A_{\varrho, r}$  in the terms of $K$--functionals of functions, generated by their radial derivatives.

%\normalsize \vskip 3mm

%\newpage
Let $L_p=L_p(\mathbb T)$, $1\le p\le \infty$, be the space of all functions $f$, given on the
%%%interval ${\mathbb T}=[0,2\pi]$,
torus ${\mathbb T}$,
with the usual norm
 $$
 \|f\|_{p}\!:=\|f\|_{L_p(\mathbb T)}
 %%%
 :=\left\{\matrix{
  \displaystyle\Big(\frac 1{2\pi}\int_{0}^{2\pi} |f(x)|dx\Big)^{1/p}\!\!,\quad \hfill &  1\le p<\infty,\cr
 \mathop{\rm ess\,sup}\limits_{x\in [0,2\pi]} |f(x)|,\, \quad \hfill &  p=\infty.}\right.
  $$

Further, let  $f\in L_1$, the Fourier coefficients of $f$ are given by
\[
\widehat f_k:=\frac 1{2\pi}\int_0^{2\pi}f(x) \mathrm{e}^{- \mathrm{i}kt}dx, \quad k\in\mathbb Z.
\]
We denote by $f\left(\varrho,x\right)$, $0\le\varrho<1$, the Poisson integral (the Poisson operator) of $f$, i.e.,
\begin{equation}\label{Poisson operator}
f\left(\varrho,x\right):=\frac{1}{2\pi}\int_0^{2\pi}f(t)P(\varrho,x-t)dt,
\end{equation}
where $P(\varrho,t)=\frac{1-\varrho^2}{|1-\varrho  \mathrm{e}^{ \mathrm{i}t}|^2}$ is the  Poisson kernel.

Leis [\ref{Leis}] considered the transformation
\[
L_{\varrho, r}(f)(x):=\sum_{k=0}^{r-1}\frac{ \mathrm{d}^k f(x)}{ \mathrm{d} n^k}\cdot\frac{(1-\varrho)^k}{k!},\quad r\in\mathbb N,
\]
where
\[
\frac{ \mathrm{d} f(x)}{ \mathrm{d} n}=\left.-\frac{\partial f(\varrho,x)}{\partial\varrho}\right|_{\varrho=1}
\]
is the normal derivative of the function $f$. He showed that if $1<p<\infty$ and
\[
\|f(\varrho, \cdot)-L_{\varrho, r}(f)(\cdot)\|_p=O\bigg(\frac{(1-\varrho)^r}{r!}\bigg),\quad\varrho\to 1-,
\]
then  $ \mathrm{d}^r/ \mathrm{d} n^r f\in L_p$.

 Butzer and Sunouchi [\ref{Butzer_Sunouchi}] considered the transformation
\[
B_{\varrho,r}(f)(x):=\sum_{k=0}^{r-1}(-1)^{\frac{k+1}{2}}f^{\{k\}}(x)\frac{(-\ln\varrho)^k}{k!},
\]
where
\[
f^{\{k\}}(x)=\left\{\matrix{
f^{(k)}, \hfill & k\in 2\mathbb Z_+,\cr
\widetilde f^{(k)}, \hfill &
%%%k+1
k-1
\in 2\mathbb Z_+.
}\right.
\]
They proved the following theorem:

{\bf Theorem A [\ref{Butzer_Sunouchi}].} {\it Assume that  $f\in L_p$, $1\le p<\infty$.

$i)$
%%%if
If
the derivatives $f^{\{j\}},$ $j=0,1,\ldots, r-1,$ are absolutely continuous  and $f^{\{r\}}\in L_p$, then
\begin{equation}\label{B-S approx}
\|f(\varrho,\cdot)-B_{\varrho, r}(f)(\cdot)\|_p=O\bigg(\frac{(-\ln\varrho)^r}{r!}\bigg),\quad\varrho\to 1-~%
%%%;
.
\end{equation}

$ii)$
%%%if
If
the derivatives $f^{\{j\}},$ $j=0,1,\ldots, r-2,$ $r\ge 2,$  are absolutely continuous, $f^{\{r-1\}}\in L_p$, $1<p<\infty,$ and relation (\ref{B-S approx}) holds, then $\widetilde f^{\{r-1\}}$ is absolutely continuous and $\widetilde f^{\{r\}}\in L_p$.
}

These results
%%%are the approximation theorem by images
summarize the approximation behaviour
of the operators $L_{\varrho, r}$ and $B_{\varrho, r}$ in the space $L_p$. %%%Better English???
In particular, Leis's result and the statement  $ii)$ of Theorem A
%%%
represent the so-called %%%are the
inverse theorems and the statement $i)$ is the
so-called %%%included
direct theorem.
Direct and inverse theorems are one of the central  theorems of approximation theory. They were studied by many authors. Here, we mention only the books [\ref{Butzer_Nessel}, \ref{DeVore Lorentz}, \ref{Trigub_Bellinsky}], which contain fundamental results in this subject.
%%% Should we include (***) already here??  (see below)
The given results are based on the investigations
%%%of
in
the papers [\ref{Butzer_Tillmann}, \ref{Butzer}], where the authors find the direct and inverse approximation theorems for the one-parameter  semi-groups of bounded
linear transformations  $\{T(t)\}$ of
some %% included
Banach space $X$ into
itself by  the ``Taylor polynomial''  %%%I changed the " "
$\sum_{k=0}^{r-1}(t^k/k!)A^kf$, where $Af$ is the infinitesimal operator of a semi--group $\{T(t)\}$.

The transformations $A_{\varrho, r}$, considered in
%%% our
this
paper, are similar to the transformations $L_{\varrho, r}$ and $B_{\varrho,r}$ as they are also  based on the ``Taylor polynomials''. %%% " "
The transformation  $A_{\varrho,r}$ are defined in the following way:

For  $\varrho\in [0,1)$, $r\in\mathbb N$ and  $f\in L_1$, we set
\begin{equation}\label{def Ar}
A_{\varrho,r}(f)(t):=\sum_{k\in {\mathbb Z}}\lambda_{|k|,r}(\varrho)\widehat f_k \mathrm{e}^{ \mathrm{i}kt},
\end{equation}
where for $k=0,1,\ldots,r-1$, the numbers $\lambda_{k,r}(\varrho)\equiv 1$  and
\begin{equation}\label{lambda for H^r}
\lambda_{k,r}(\varrho):=\sum_{j=0}^{r-1}
{k\choose j}(1-\varrho)^j\varrho^{k-j}, k=r,r+1,\ldots,\quad\varrho\in[0,1].
\end{equation}

The transformation  $A_{\varrho,r}$ can be considered as
%%% the
a
linear operator
%%%of
on
 $L_1$ into
itself. Indeed, $\lambda_{k,r}(0){=}0$ and for all  $k=r,r+1,\ldots$ and $\varrho\in(0,1)$, we have
\[
\sum_{j=0}^{r-1}
{k\choose j}(1-\varrho)^j\varrho^{k-j}\le
rq^{k}k^{r-1},~\mbox{where}~0<q:=\max\{1-\varrho,\varrho\}<1.
\]
Therefore, for any function $f\in L_1$ and for any $0<\varrho<1$, the series on the right-hand side of (\ref{def Ar}) is majorized by the convergent series $2r\|f\|_1\sum_{k=r}^{\infty}q^{k}k^{r-1}$.

Note that if the function  $f\in L_1$ and it has the Fourier series of power type, i.e.,
$
f(x)\sim\sum_{k=0}^\infty\widehat f_k\mathrm{e}^{ \mathrm{i}kx},
$
then
$
f(\varrho, x)=f(z):=\sum_{k=0}^\infty\widehat f_kz^k,\quad z=\varrho \mathrm{e}^{ \mathrm{i}x}.
$

The
%%%relations of the operators $A_{\varrho, r}$ and  the "Taylor polynomials"\  are
relation between the operators $A_{\varrho, r}$ and  the ``Taylor polynomials''  is
 shown in the following statement.

{\bf Lemma 1.} {\it   Assume that  $f\in L_1$. Then for any numbers $r\in \mathbb N,$ $\varrho\in[0,1)$ and $x\in
%%%[0,2\pi]
\mathbb T$,
\begin{equation}\label{A_P}
A_{\varrho, r}(f)(x)=\sum_{k=0}^{r-1}\frac{\partial^k f\left(\varrho,x\right)}{\partial\varrho^k}\cdot\frac{(1-\varrho)^k}{k!}.
\end{equation}
}

{\bf Proof.} Let us associate the function $f$ with the functions
\begin{equation}\label{NEW1}
f_1(z):=\widehat f_0/2+\sum_{k=1}^\infty\widehat f_kz^k\quad {\rm and}\quad f_2(z):=\widehat f_0/2+\sum_{k=1}^\infty\widehat f_{-k}z^k,
 \end{equation}
which are holomorphic in the disc $\mathbb D:=\{z\in {\mathbb C}:|z|< 1\}$.

From Lemma 4
in %%%added
[\ref{Savchuk}], it follows that for any $z\in\overline{\mathbb D}$,
\begin{equation}\label{NEW2}
\frac{\widehat f_0}{2}+\sum_{k=1}^{r-1}\widehat f_kz^k+\sum_{k=r}^\infty\lambda_{k,r}(\varrho)\widehat f_kz^k=\frac{\widehat f_0}{2}+\sum_{k=1}^{r-1}z^kf^{(k)}_1(\varrho z)\frac{(1-\varrho)^k}{k!}
 \end{equation}
and
\begin{equation}\label{NEW3}
\frac{\widehat f_0}{2}+\sum_{k=1}^{r-1}\widehat f_{-k}\overline z^k+\sum_{k=r}^\infty\lambda_{k,r}(\varrho)\widehat f_{-k}\overline z^k=\frac{\widehat f_0}{2}+\sum_{k=1}^{r-1}\overline z^kf^{(k)}_2(\varrho \overline z)\frac{(1-\varrho)^k}{k!},
 \end{equation}
where for $r=1$, the sums $\sum_{k=1}^0$ are assumed to be zero.

Actually, in [\ref{Savchuk}],  the relations of the kind as in (\ref{NEW2}) and (\ref{NEW3}) were proved for $z\in {\mathbb D}$, but  such restrictions are not important.

Adding these two equalities at $z=\mathrm{e}^{ \mathrm{i}x}$  and taking into account the relation
\begin{equation}\label{f1+f2=P}
\mathrm{e}^{ \mathrm{i}kx}f^{(k)}_1(\varrho \mathrm{e}^{ \mathrm{i}x})+\mathrm{e}^{-ikx}f^{(k)}_2(\varrho \mathrm{e}^{-ix})=\frac{\partial^k f\left(\varrho,x\right)}{\partial\varrho^k},
\end{equation}
we get (\ref{A_P})%
%%%. Lemma is proved.
, which proves the Lemma.

%%%%%%%%%%%%%%%%%%%%%%%%%%%%%%%%%%%%%%%%%%%%%%%%%%%%%%%%%%%%%%%%%%%%%%%%%%%%%%%%%%%%%%%%%%%%%%%%%%%%%%%%%%%%%%%%

Now let us formulate direct and inverse approximation theorems by the operators $A_{\varrho, r}$  in the terms of $K$--functionals of functions, generated by their radial derivatives.
%%% Should we write the following two sentences already earlier?? I propose the position marked with (***)  Ok

%%%

 Let us give all necessary definitions. If for a function $f\in L_1$ and for a positive integer  $n$, there exists the function $g\in L_1$  such that
\[
\widehat g_k=\left\{\matrix{ 0,\hfill  & \mbox{\rm if}\quad |k|<n, \cr  {\displaystyle \frac{|k|!}{(|k|-n)!}}\widehat f_k,\hfill & \mbox{\rm if}\quad |k|\ge n,}\right.,\quad  k\in\mathbb Z
\]

then we say that for the function $f$, there exists the radial derivative  $g$ of order $n$, for which we use the notation $f^{[n]}$. Here, we use the term ``radial derivative'' %%" "
in view of the following fact.

If the function $f^{[r]}\in L_1$, then its Poisson integral can be presented
%%%in the form
as
\begin{equation}\label{diff f[n]}
f^{[r]}(\varrho,x)=(f(\varrho,\cdot))^{[r]}(x)=\varrho^r\frac{\partial^r f\left(\varrho,x\right)}{\partial\varrho^r}
\quad \varrho\in[0,1),~\forall~x\in
%%%[0,2\pi]
\mathbb T.
\end{equation}
Hence, by virtue of the theorem of  limit values of Poisson  integral (see, for example, [\ref{Rudin}, p.~27]), for almost all $x\in
%%%[0, 2\pi]
\mathbb T
$, we have $f^{[r]}(x)=\lim_{\varrho\to 1-}f^{[r]}(\varrho,x)$.

Relation (\ref{diff f[n]}) can be
%%%easy
easily
proved by term by term differentiation with respect to the variable $\varrho$ of the  decomposition of Poisson integral into the uniformly convergent series
\begin{equation}\label{series fo Poisson}
f\left(\varrho,x\right)=\sum_{k\in\mathbb Z}\varrho^{|k|}\widehat f_k\mathrm{e}^{ \mathrm{i}kx}\quad\forall~\varrho\in[0,1),~x\in
%%%[0, 2\pi]
\mathbb T.
\end{equation}

From the definition of radial derivative, in particular, it follows the
%%%following
differentiation rule:

\noindent{\it
%%%if
If $f(x)=\sum_{|k|\le m}\widehat f_k\mathrm{e}^{ \mathrm{i}kx}$, $m\in\mathbb Z_+,$ then
 \begin{equation}\label{rule for differ}
 f^{[n]}(x)=\left\{\matrix{
 0,\hfill  & \mbox{\rm if}\quad m<n,\cr
 \displaystyle \sum_{n\le |k|\le m}\frac{|k|!}{(|k|-n)!}\widehat f_k\mathrm{e}^{ \mathrm{i}kx},\hfill   & \mbox{\rm if}\quad m\ge n.
}\right.
 \end{equation}}

 %%%%%%%%%%%%%%%%%%%%%%%%%%%%%%%%%%%%%%%%%%%%%%%%%%%%%%%%%%%%%%%%%%%%%%%%%%%%%%%%%%%%%%%%%%%%%%%%%%%%%%%%%%%

% and \[ f^{[r]}(\varrho, x)=f^{(r)}(z),\quad z=\varrho \mathrm{e}^{ \mathrm{i}x}. \]

%%%%%%%%%%%%%%%%%%%%%%%%%%%%%%%%%%%%%%%%%%%%%%%%%%%%%%%%%%%%%%%%%%%%%%%%%%%%%%%%%%%%%%%%%%%%%%%%%%%%%%%%%%%

In the space $L_p$,
the %%% added
$K$--functional of
%%%the
a
function $f$ (see, for example, [\ref{DeVore Lorentz}, Chap.~6])  generated by the radial derivative of order  $n$,  is the following quantity:
\[
K_n(\delta, f)_p:=\inf\left\{\left\|f-h\right\|_p+\delta^n\left\|h^{[n]}\right\|_p: h^{[n]}\in L_p\right\},\quad\delta>0.
\]

Further, we consider the functions $\omega(t)$, $t\in [0,1]$, satisfying the following conditions:

1) $\omega(t)$ is continuous on $[0,1]$;

2) $\omega(t)\uparrow$;

3) $\omega(t)\not=0$ for any $t\in (0,1]$;

4) $\omega(t)\to 0$ as $t\to 0$;

\noindent
%%% as
and the %%% I think it is now the right English for that you wanted to say??
well-known Zygmund--Bari--Stechkin conditions (see, for example, [\ref{Bari_Stechkin}]):
$$%\begin{equation}\label{Bari-Steckin condition-1}
({\mathscr Z}) \quad\quad\quad\quad\int_0^\delta\frac{\omega(t)}{t}dt=O(\omega(\delta)),\quad\quad\delta>0,\quad\quad\quad
$$%\end{equation}
$$%\begin{equation}\label{Bari-Steckin condition-2}
({\mathscr Z}_n) \quad\quad\quad \int_\delta^1\frac{\omega(t)}{t^{n+1}}dt=O\bigg(\frac{\omega(\delta)}{\delta^n}\bigg),\quad\delta>0,\ n\in {\mathbb N}.
$$%\end{equation}

%%% Should we refer from now now on to (13) and (14)  or better to (Z) and (Z_n)??  Please change if you agree.

The main results of
%%%the paper are contained in the following statements:
this paper are contained in the following two statements:

{\bf Theorem 1.} {\it Assume that $f\in L_p,$ $1\le p\le\infty$, $n, r\in\mathbb N$, $n\le r$ and the function  $\omega(t)$, $t\in [0,1]$, satisfies conditions 1)--4) and $({\mathscr Z})$ . If
\begin{equation}\label{K-funct est}
K_{n}\left(\delta, f^{[r-n]}\right)_p=O(\omega(\delta)),\quad\delta\to 0+,
\end{equation}
then
\begin{equation}\label{f-Ap est}
\|f-A_{\varrho, r}(f)\|_p=O\left((1-\varrho)^{r-n}\omega(1-\varrho)\right),\quad\varrho\to 1-.
\end{equation}
}

{\bf Theorem 2.} {\it  Assume that $f\in L_p,$ $1\le p\le\infty$, $n, r\in\mathbb N$, $n\le r$ and the function  $\omega(t)$, $t\in [0,1]$, satisfies conditions 1)--4), $({\mathscr Z})$  and $({\mathscr Z}_n)$ . If relation (\ref{f-Ap est}) holds, then $f^{[r-n]}\in L_p$ and  relation (\ref{K-funct est}) also  holds.
}

We note that in the case where $\omega(t)$ is a power function: $\omega(t)=t^\alpha$, $\alpha>0$, the results of the Theorems 1 and 2  were announced in [\ref{Prestin_Savchuk_Shidlich}].

{\bf\textsl {Remark} 1.} For a given $n\in\mathbb N$, from condition $({\mathscr Z}_n)$  it follows that  $\mathop{\rm lim~inf}\limits_{\delta\to 0+}(\delta^{-n}\omega(\delta))>0$ or, equivalently, that $(1-\varrho)^{r-n}\omega(1-\varrho)\ll (1-\varrho)^r$ as $\varrho\to~1-$. Therefore, if condition $({\mathscr Z}_n)$  is satisfied, then the quantity on the right-hand side of (\ref{f-Ap est}) decreases to zero as $\varrho\to 1-$ not faster, than the function $(1-\varrho)^r$. Also note that the relation $
\|f-A_{\varrho, r}(f)\|_p=o\left((1-\varrho)^r)\right),\ \varrho\to 1-,
$
 holds only in the trivial case when $f(x)=\sum_{|k|\le r-1}\widehat f_k\mathrm{e}^{ \mathrm{i}kx}$, and in such case, the theorems are
 easily %%%added
 true. This fact is   related to the so-called  saturation property of the approximation method,  generated by the operator $A_{\varrho, r}$. In particular, in [\ref{Savchuk}], it was shown that the operator $A_{\varrho, r}$  generates the linear approximation method of holomorphic functions, which is saturated in the space $H_p$ with the saturation order $(1-\varrho)^r$ and the saturation class $H^{r-1}_p\mathop{\rm Lip}1$.

Before proving
%%%of the theorems
the Theorems
1 and 2, let us give some auxiliary results.
For any   $f\in L_1$, $1\le p\,{\le}\infty$, $0\le\varrho<1$ and  $r=0,1,2,\ldots$, we set
\begin{equation}\label{M_p_diff}
M_p(\varrho, f,r):=\varrho^r\Big\| \frac{\partial^r f\left(\varrho,\cdot\right)}{\partial\varrho^r}\Big\|_{p}=
\Big\|(f(\varrho,\cdot))^{[r]}(\cdot)\Big\|_{p}.
\end{equation}

{\bf Lemma 2.} {\it Assume that  $f\in L_p$, $1\le p\le \infty$. Then for any numbers $n\in \mathbb N$ and $\varrho\in[1/2,1)$,
%%% change of equation style
\begin{eqnarray*}
\frac{1}{2n!}(1-\varrho)^{n}M_p\left(\varrho, f,n\right)&\le &K_{n}\left(1-\varrho, f\right)_p
\\
&\le&
\|f-A_{\varrho, n}(f)\|_p+\frac{4^{n}-1}{3}(1-\varrho)^{n}M_p\left(\sqrt{\varrho}, f,n\right).
\end{eqnarray*}}

{\bf Proof.} First, let us note that the statement of Lemma 2 is trivial in the case,
%%%when
if
$f$ is a trigonometric polynomial of order not exceeding $n-1$, i.e.,
%%%when
if
$f(x)=\sum_{|k|\le n-1}\widehat f_k \mathrm{e}^{ \mathrm{i}kx}$, as well as in the case,
%%%when
if $\varrho=0$. Therefore, further in the proof, we exclude these two cases.

Let $g$ be a function such that  $g^{[n]}\in L_p$.

Since
\[
\frac{1-\varrho^2}{|1-\mathrm{e}^{ \mathrm{i}(x-t)}\varrho|^2}=\frac{1}{1-\mathrm{e}^{ \mathrm{i}(x-t)}\varrho}+\frac{1}{1-e^{-i(x-t)}\varrho}-1,
\]
then by virtue of  (\ref{Poisson operator}), for any numbers $\varrho\in[0,1)$ and $x\in
%%%[0, 2\pi]
\mathbb T$,
we have
\[%\begin{eqnarray*} %%% new equationstyle
\frac{\partial^n f\left(\varrho,x\right)}{\partial\varrho^n}
=\frac{1}{2\pi}\int_0^{2\pi}\left(f(t)-g(t)\right)
\frac{\partial^n }{\partial\varrho^n}\left(\frac{1-\varrho^2}{|1-\mathrm{e}^{ \mathrm{i}(x-t)}\varrho|^2}\right)dt+
\frac{\partial^n g(\varrho,x)}{\partial\varrho^n}
\]\[
=\frac{n!}{2\pi}\int_0^{2\pi}\!\!\!\left(f(t)-g(t)\right)\left(\frac{\mathrm{e}^{ \mathrm{i}r(x-t)}}{(1-\mathrm{e}^{ \mathrm{i}(x-t)}\varrho)^{n+1}}+
\frac{\mathrm{e}^{- \mathrm{i}r(x-t)}}{(1-\mathrm{e}^{- \mathrm{i}(x-t)}\varrho)^{n+1}}\right)dt+
\frac{\partial^n g(\varrho,x)}{\partial\varrho^n}
\]\[
=\frac{n!}{\pi}\int_0^{2\pi}\left(f(t)-g(t)\right)\mathop{\rm Re}\frac{\mathrm{e}^{ \mathrm{i}r(x-t)}}{(1-\mathrm{e}^{ \mathrm{i}(x-t)}\varrho)^{n+1}}dt+
\frac{\partial^n g(\varrho,x)}{\partial\varrho^n}.
\]%\end{eqnarray*}
Hence, changing the variables of integration and using the integral Minkowski inequality, we obtain
\begin{eqnarray*}
\bigg\|\frac{\partial^n f\left(\varrho,\cdot\right)}{\partial\varrho^n}\bigg\|_p&\le&\frac{n!}{\pi}\int_0^{2\pi}\frac{dt}{|1-\varrho \mathrm{e}^{ \mathrm{i}t}|^{n+1}}\|f-g\|_p+\bigg\|\frac{\partial^n g\left(\varrho,\cdot\right)}{\partial\varrho^n}\bigg\|_p\\ &\le&
\frac{2n!}{(1-\varrho)^n}\|f-g\|_p+
\bigg\|\frac{\partial^n g\left(\varrho,\cdot\right)}{\partial\varrho^n}\bigg\|_p.
\end{eqnarray*}

Taking into account (\ref{diff f[n]}), (\ref{M_p_diff})  and the inequality $\|g^{[n]}(\varrho,\cdot)\|_p\le\|g^{[n]}\|_p$, we see that for any $\varrho\in (0,1)$,
\[
\frac{1}{2n!}(1-\varrho)^n M_p\left(\varrho, f,n\right)\le\left\|f-g\right\|_p+(1-\varrho)^n\left\|g^{[n]}\right\|_p.
\]
Considering the infimum over all functions $g$ such that $g^{[n]}\in L_p$, we conclude that
\[
\frac{1}{2n!}(1-\varrho)^nM_p\left(\varrho, f,n\right)\le K_{n}\left(1-\varrho, f\right)_p.
\]

On the other hand, from the definition of the  $K$--functional, it follows that
\begin{equation}\label{K_n}
K_n\left(1-\varrho, f\right)_p\le\|f-A_{\varrho, n}(f)\|_p+(1-\varrho)^n\left\|\left(A_{\varrho, n}(f)\right)^{[n]}\right\|_p.
\end{equation}
According to (\ref{A_P}) and (\ref{diff f[n]}), we have
\begin{eqnarray*}
(A_{\varrho,n}(f))^{[n]}(x)&=&\bigg(\sum_{k=0}^{n-1}\frac{(f(\varrho,\cdot))^{[k]}(\cdot)}{\varrho^k k!}(1-\varrho)^k\bigg)^{[n]}(x)
\\ &=&\sum_{k=0}^{n-1}\frac{((f(\varrho,\cdot))^{[k]}(\cdot))^{[n]}(x)}{\varrho^k k!}(1-\varrho)^k.
\end{eqnarray*}

Since for any nonnegative integers $k$ and $n$
\begin{equation}\label{diff f[n]_f[k]}
((f(\varrho,\cdot))^{[n]}(\cdot))^{[k]}(x)=((f(\varrho,\cdot))^{[k]}(\cdot))^{[n]}(x),
\end{equation}
%%%then
we obtain
$$
(A_{\varrho, n}(f))^{[n]}(x)=\sum_{k=0}^{n-1}\frac{((f(\varrho,\cdot))^{[n]}(\cdot))^{[k]}(x)}{\varrho^k k!}(1-\varrho)^k.
$$
This yields
\begin{equation}\label{A11}
\|(A_{\varrho, n}(f))^{[n]}\|_p\le \sum_{k=0}^{r-1}\frac{\|((f(\varrho,\cdot))^{[n]}(\cdot))^{[k]}\|_p}{\varrho^k k!}(1-\varrho)^k.
\end{equation}

By virtue of the definition of the Poisson integral, for any $k=0,1,\ldots,r-1$, we have
\begin{eqnarray*}
((f(\varrho,\cdot))^{[n]}(\cdot))^{[k]}(x)&=&\bigg(\sum\limits_{|j|\ge n}\frac{|j|!}{(|j|-n)!}\widehat{f}_j\varrho^\frac{|j|}2\mathrm{e}^{ \mathrm{i}jx}\varrho^\frac{|j|}2\bigg)^{[k]}(x) \\
&=&\bigg(\frac 1{2\pi}\int\limits_0^{2\pi} (f(\sqrt{\varrho},\cdot))^{[n]}(t)P(\sqrt{\varrho}, t-\cdot)dt\bigg)^{[k]}(x)\\ &=&\frac 1{2\pi}\int\limits_0^{2\pi} (f(\sqrt{\varrho},\cdot))^{[n]}(t)\sum\limits_{|\nu|\ge k}\frac{|\nu|!}{(|\nu|-k)!}\varrho^\frac{|\nu|}2\mathrm{e}^{ \mathrm{i}\nu(t-x)}dt\\
&=&\frac 1{2\pi}\! \int\limits_0^{2\pi}\! \!  (f(\sqrt{\varrho},\cdot))^{[n]}(t+x)\! \sum\limits_{|\nu|\ge k}\frac{|\nu|!}{(|\nu|-k)!}\varrho^\frac{|\nu|}2\mathrm{e}^{ \mathrm{i}\nu t}dt\\ &=&
\frac 1{2\pi}\! \int\limits_0^{2\pi} \! \! (f(\sqrt{\varrho},\cdot))^{[n]}(t+x)\bigg(\tau^k  \frac{\partial^k }{\partial\tau^k}P(\tau,t) \bigg) \bigg|_{\tau=\sqrt{\varrho}}\! \! dt.
\end{eqnarray*}
Using  the integral Minkowski inequality, for  $k=0$, we obtain
\begin{equation}\label{A12}
\|((f(\varrho,\cdot))^{[n]}(\cdot))^{[k]}\|_p=\|(f(\varrho,\cdot))^{[n]}\|_p
$$
$$
\le
M_p(\sqrt{\varrho},f,n)\frac 1{2\pi}\int\limits_0^{2\pi}|P(\sqrt{\varrho},t)|dt=M_p(\sqrt{\varrho},f,n).
\end{equation}
If $k=1,2,\ldots$, then
 $$
 \frac{\partial^k }{\partial\tau^k}P(\tau,t)=  \frac{\partial^k }{\partial\tau^k}\bigg(\frac 1{1-\tau \mathrm{e}^{ \mathrm{i}t}}+\frac {\tau \mathrm{e}^{- \mathrm{i}t}}{1-\tau \mathrm{e}^{- \mathrm{i}t}}\bigg)=\frac{k!\,\mathrm{e}^{ \mathrm{i}kt}}{(1-\tau \mathrm{e}^{ \mathrm{i}t})^{k+1}}+\frac{k!\,\mathrm{e}^{- \mathrm{i}kt}}{(1-\tau \mathrm{e}^{- \mathrm{i}t})^{k+1}}.
 $$
This similarly yields
\begin{equation}\label{A13}
\|(f^{[n]}(\varrho, \cdot))^{[k]}\|_p \le M_p(\sqrt{\varrho},f,n)\frac 1{2\pi}\int\limits_0^{2\pi}\bigg|\bigg(\tau^k  \frac{\partial^k }{\partial\tau^k} P(\tau,t)\bigg) \bigg|_{\tau=\sqrt{\varrho}}\bigg|dt
$$
$$
\le  2k!\, M_p(\sqrt{\varrho},f,n)\frac 1{2\pi}\int\limits_{0}^{2\pi}
\frac {dt}{|1-\sqrt{\varrho}\mathrm{e}^{ \mathrm{i}t}|^{k+1}}\le \, M_p(\sqrt{\varrho},f,n)\frac {2^{k}k!}{(1-\varrho)^{k}}. %%%end of sentence!
\end{equation}
Combining relations (\ref{A11})--(\ref{A13}), we see that for any $\varrho\in[1/2,1)$,
\begin{eqnarray}\nonumber
 \int_\varrho^1 \|(A_{\varrho, n}(f))^{[n]}\|_p &\le& M_p(\sqrt{\varrho},f,n)+M_p(\sqrt{\varrho},f,n)\sum_{k=1}^{n-1}4^k\\  \label{A_{varrho, n}}
&=& M_p(\sqrt{\varrho},f,n)\frac{4^{n}-1}3.
\end{eqnarray}

By virtue
of %%% added
(\ref{A_{varrho, n}}) and  (\ref{K_n}), we conclude that
$$
K_n\left(1-\varrho, f\right)_p\le\|f-A_{\varrho, n}(f)\|_p+\frac{4^{n}-1}3(1-\varrho)^nM_p(\sqrt{\varrho},f,n) %%%.
$$
%%%Lemma is proved.
which proves the Lemma.

{\bf Lemma 3.}  {\it Assume that $r\in {\mathbb N}$, $1\le p\le \infty$  and $\varrho\in[1/2,1)$. Then for any function $f\in L_p$
\begin{equation}\label{|A_r(f)|}
\|(A_{\varrho, r}(f))^{[r]}\|_p \le C_r\frac{\|f\|_p}{(1-\varrho)^r},
\end{equation}
where the
%%%quantity
constant
$C_r$ depends only on $r$.
}

{\bf Proof.} By virtue of (\ref{diff f[n]}),  for any function $f\in L_p$ and all $x\in
%%%[0, 2\pi]
\mathbb T$, we have
\begin{eqnarray*}
(f(\varrho,\cdot))^{[r]}(x)
&=&\frac{\varrho^r}{2\pi}\int_0^{2\pi}f(t)\frac{\partial^r}{\partial\varrho^r}\left(\frac{1-\varrho^2}{|1-\mathrm{e}^{ \mathrm{i}(x-t)}
\varrho|^2}\right)dt\\
&=&\frac{r!\varrho^r}{2\pi}\int_0^{2\pi}f(t)\left(\frac{\mathrm{e}^{ \mathrm{i}r(x-t)}}{(1-\mathrm{e}^{ \mathrm{i}(x-t)}\varrho)^{r+1}}+
\frac{\mathrm{e}^{- \mathrm{i}r(x-t)}}{(1-\mathrm{e}^{- \mathrm{i}(x-t)}\varrho)^{r+1}}\right)dt\\
&=&\frac{r!\varrho^r}{\pi}\int_0^{2\pi}f(t)\mathop{\rm Re}\frac{\mathrm{e}^{ \mathrm{i}r(x-t)}}{(1-\mathrm{e}^{ \mathrm{i}(x-t)}\varrho)^{r+1}}dt.
\end{eqnarray*}
Making a change of variables of integration and using the integral Minkowski inequality, we obtain
\begin{equation}\label{|f^[r]|}
M_p\left(\varrho, f,r\right)\le\frac{r!}{\pi}\int_0^{2\pi}\frac{dt}{|1-\varrho \mathrm{e}^{ \mathrm{i}t}|^{r+1}}\|f\|_p\le\frac{2r!}{(1-\varrho)^r}\|f\|_p.
\end{equation}
Combining this relation and relation (\ref{A_{varrho, n}})  with $n=r$, we conclude that
\[
\|(A_{\varrho, r}(f))^{[r]}\|_p\le M_p(\sqrt{\varrho},f,r)\frac{4^{r}-1}3\le \frac{2r!(4^{r}-1)}{3(1-\sqrt{\varrho})^r}\|f\|_p
$$
$$
\le \frac {r!(2^{3r+1}-2^{r+1})}{3}\cdot \frac{\|f\|_p}{(1-\varrho)^r}.
\]

{\bf Lemma 4.} {\it Assume that $r\in {\mathbb N}$ and $0\le \varrho<1$. Then for any function $f\in L_p$, $1\le p\le \infty$, such that
\begin{equation}\label{f-A(f)Condition}
\int_\varrho^1 \bigg\|\frac{\partial^rf(\zeta,\cdot)}{\partial\zeta^r}\bigg\|_p(1-\zeta)^{r-1}d\zeta<\infty
\end{equation}
and for almost all $x\in
%%%[-\pi, \pi]
\mathbb T$,
\begin{equation}\label{f-A(f)}
f(x)-A_{\varrho, r}(f)(x)=\frac{1}{(r-1)!}\int_\varrho^1 \frac{\partial^rf(\zeta,x)}{\partial\zeta^r}(1-\zeta)^{r-1}d\zeta.
\end{equation}

}

{\bf Proof.}  For  fixed $r\in {\mathbb N}$ and $0\le \varrho<1$, the integral on the right-hand side of (\ref{f-A(f)}), defines a certain function $F(x)$. By virtue of (\ref{f-A(f)Condition}) and  the integral Minkowski inequality, we conclude that the function  $F$ belongs to the space $L_p$. Let us find the Fourier coefficients of $F$ and compare them with
 the Fourier coefficients of the function $G:=f-A_{\varrho, r}(f)$. Since for $r\in {\mathbb N}$, $$
\frac{\partial^rf(\zeta,x)}{\partial\zeta^r}=\sum\limits_{|k|\ge r} \frac{|k|!}{(|k|-r)!}\widehat{f}_k \zeta^{|k|-r}\mathrm{e}^{ \mathrm{i}kx},
$$
then $\widehat{F}_{k}=0$, when $|k|<r$. If $|k|\ge r$, then integrating by parts, we see that
\begin{equation}\label{Fouries coeff}
\widehat{F}_{k}=\frac 1{2\pi}\int_0^{2\pi}F(t)\mathrm{e}^{- \mathrm{i}kt}dt=\widehat{f}_k \sum\limits_{j=r}^{|k|} {|k|\choose j}(1-\varrho)^j\varrho^{|k|-j}. %%% end of sentence!
\end{equation}

On the other hand, if $|k|<r$
the Fourier coefficients $\widehat{G}_{k}$ of the function $G$ are equal to zero.
%%%when $|k|<r$.
If $|k|\ge r$, then in view of the equality
 $$
 \sum_{j=0}^{k}{k\choose j}(1-\varrho)^j\varrho^{k-j}=\big((1-\varrho)+\varrho\big)^{k}=1, \qquad  k=0,1,\ldots,
 $$
we see that
 $$
 \widehat{G}_{k}=(1-\lambda_{|k|,r}(\varrho))\widehat{f}_k=\widehat{f}_k \sum\limits_{j=r}^{|k|} {|k|\choose j}(1-\varrho)^j\varrho^{|k|-j}.
 $$
Therefore, for all $k\in {\mathbb Z}$, we have $\widehat{F}_{k}=\widehat{G}_{k}$. Hence, for almost all $x\in
%%%[-\pi, \pi]
\mathbb T$, relation (\ref{f-A(f)}) holds.

%%%%%%%%%%%%%%%%%%%%%%%%%%%%%%%%%%%%%%%%%%%%%%%%%%%%%%%%%%%%%%%%%%%%%%%%%%%%%%%%%%%%%%%%%%%%%%%%%%%%%%%%%%%%%%%%%%%%%%%%%%%%%%%%%%%%%%%%%

\bigskip
{\bf Proof of Theorem 1.} Assume that the function $f$ is such that $f^{[r-n]}\in L_p$ and relation (\ref{K-funct est}) is satisfied.  Let us apply  the first inequality
%%%in
of
Lemma 2 to the function $f^{[r-n]}$. In view of  (\ref{diff f[n]}) and (\ref{M_p_diff}), we obtain
\[
\frac{1}{2n!}(1-\varrho)^{n}M_p\left(\varrho, f,r\right)\le K_{n}\left(1-\varrho, f^{[r-n]}\right)_p.
\]
This yields
\begin{equation}\label{est Mp}
M_p\left(\varrho, f,r\right)\le C\frac{\omega(1-\varrho)}{(1-\varrho)^n}, \quad \varrho\to 1-.
\end{equation}

Using  relations (\ref{M_p_diff}), (\ref{est Mp}) and $({\mathscr Z})$  and the integral Minkowski inequality, we obtain
 \begin{eqnarray}\nonumber
 \int_\varrho^1 \bigg\|\frac{\partial^rf(\zeta,\cdot)}{\partial\zeta^r}\bigg\|_p(1-\zeta)^{r-1}d\zeta&\le&\int_\varrho^1 M_p\left(\zeta, f,r\right)\frac{(1-\zeta)^{r-1}}{\zeta^r}d\zeta\\ \nonumber
&\le& {2^{r}C}
(1-\varrho)^{r-n}\int_\varrho^1\frac{\omega(1-\zeta)}{1-\zeta}d\zeta\\  \label{est Mp estim}
&=& \!\!O\left((1-\varrho)^{r-n}\omega(1-\varrho)\right), \ \varrho\to 1-.
\end{eqnarray}
Therefore, for almost all $x\in
%%%[-\pi, \pi]
\mathbb T$, relation (\ref{f-A(f)}) holds. Hence, by virtue of  (\ref{f-A(f)}),
using the integral Minkowski inequality and (\ref{est Mp estim}),  we
finally %%% added
get (\ref{f-Ap est}):
 \begin{eqnarray*}
\|f-A_{\varrho, r}(f)\|_p&\le&\frac{1}{(r-1)!}\int_\varrho^1 M_p\left(\zeta, f,r\right)\frac{(1-\zeta)^{r-1}}{\zeta^r}d\zeta \\ &=&O\left((1-\varrho)^{r-n}\omega(1-\varrho)\right), \quad \varrho\to 1-.
\end{eqnarray*}

\bigskip

{\bf Proof of Theorem 2.} First, let us note that for any function $f\in L_p$, $1\le p\le\infty$, and all fixed numbers $s,r\in {\mathbb N}$ and $\varrho\in (0,1)$
 \begin{eqnarray*}
\|A_{\varrho, r}^{[s]}(f)\|_p&=&\Big\|\sum_{|k|\ge s} \frac{|k|!}{(|k|-s)!}\omega_{|k|}(\varrho)\widehat f_k\mathrm{e}^{ \mathrm{i}kt}\Big\|_p \\ &\le&2r\|f\|_p\bigg(C+\sum_{k\ge \max\{s,r\}}q^{k}k^{s+r-1}\bigg)<\infty,
\end{eqnarray*}
where $0<q=\max\{1-\varrho,\varrho\}<1$.

Put $\varrho_k:=1-2^{-k},~k\in{\mathbb N},$ and $A_k:=A_k(f):=A_{\varrho_k,r}(f)$. For any $x\in
%%%[0, 2\pi]
\mathbb T$
and $s\in {\mathbb N}$, consider the series
\begin{equation}\label{series1}
A_0^{[s]}(f)(x)+\sum\limits_{k=1}^\infty (A_k^{[s]}(f)(x)-A_{k-1}^{[s]}(f)(x)).
\end{equation}
According to the definition of the operator  $A_{\varrho,r}$, we see that for any $\varrho_1, \varrho_2\in[0,1)$ and $r\in {\mathbb N}$,
\[
A_{\varrho_1,r}\left(A_{\varrho_2,r}(f)\right)=A_{\varrho_2,r}\left(A_{\varrho_1,r}(f)\right).
\]
By virtue of Lemma 3 and relation (\ref{f-Ap est}), for any  $k\in {\mathbb N}$ and  $s\in {\mathbb N}$, we have
\begin{equation}\label{OA}
\left\|A^{[s]}_{k}-A^{[s]}_{k-1}\right\|_p
=\left\|A^{[s]}_{k}(f-A_{k-1}(f))-A^{[s]}_{k-1}(f-A_{k}(f))
\right\|_p
$$
$$
\le\left\|A^{[s]}_{k}(f-A_{k-1}(f))\right\|_p+\left\|A^{[s]}_{k-1}(f-A_{k}(f)) \right\|_p
$$
$$
\le C_s\frac{\left\|f-A_{k-1}(f)\right\|_p}{(1-\varrho_k)^{s}}+C_s\frac{\left\|f-A_{k}(f)\right\|_p}
{(1-\varrho_{k-1})^{s}}
 $$
 $$
 =
O\left(\frac{\omega(1-\varrho_{k-1})}{(1-\varrho_{k})^{s-r+n}}\right)+
O\left(\frac{\omega(1-\varrho_{k})}{(1-\varrho_{k-1})^{s-r+n}}\right),\quad k\to +\infty.
\end{equation}
Therefore, for any $s\le r-n$,
\begin{equation}\label{series12}
\left\|A^{[s]}_{k}-A^{[s]}_{k-1}\right\|_p=O\left(\omega(1-\varrho_{k-1})\right)=
O\left(\omega(2^{-(k-1)})\right),\quad k\to +\infty.
\end{equation}
Consider the sum  $\sum_{k=1}^N \omega(2^{-(k-1)})$, $N\in {\mathbb N}$. Taking into account the monotonicity of the function $\omega$ and $({\mathscr Z})$ , we see that for all $N\in {\mathbb N}$,
\begin{eqnarray}\nonumber
 \sum_{k=1}^N\omega(2^{-(k-1)})&\le& \omega(1)+\int_1^{N} \omega(2^{-(t-1)})dt
\\ \label{summe1}
&=& \omega(1)+ \frac 1{\ln 2}\int_{2^{-N+1}}^1 \frac{\omega(\tau)}{\tau} d\tau\le C\omega(1)<\infty.
\end{eqnarray}

Combining relations (\ref{series12}) and (\ref{summe1}), we  conclude that for all $1\le p\le\infty$, the series in (\ref{series1}) converges in the
%%%metrics
norm
of the space $L_p$. Hence, by virtue of the Banach--Alaoglu theorem, for any  $s=0,1,\ldots,  r-n$, there exists the subsequence
\begin{equation}\label{series123}
S^{[s]}_{N_j}(x)=A_0^{[s]}(f)(x)+\sum\limits_{k=1}^{N_j} (A_k^{[s]}(f)(x)-A_{k-1}^{[s]}(f)(x)),\quad j=1,2,\ldots
\end{equation}
of partial sums of this series, converging to a certain function  $g\in L_p$ almost everywhere on
%%%$[0, 2\pi]$
$\mathbb T$
 as $j\to\infty$.

Let us show that $g=f^{[s]}$. For this, let us find the Fourier coefficients of the function $g$. For any fixed  $k\in\mathbb Z$ and all $j=1,2,\ldots,$ we have
\[
\widehat g_k:=\frac 1{2\pi}\int_0^{2\pi}S^{[s]}_{N_j}(t)\mathrm{e}^{- \mathrm{i}kt}{dt}+
\frac 1{2\pi}\int_0^{2\pi}(g(t)-S^{[s]}_{N_j}(t))\mathrm{e}^{- \mathrm{i}kt}{dt}.
\]
Since the sequence $\{S^{[s]}_{N_j}\}_{j=1}^\infty$ converges almost everywhere on
%%%$[0, 2\pi]$
$\mathbb T$
to the function $g$, then the second integral on the right-hand side of the last equality tends to zero as $j\to\infty$. By virtue of (\ref{series123}) and the definition of the radial derivative,
%%%the first integral is equal to zero, when $|k|<s$,
for $|k|<s$ the first integral is equal to zero,
and for all $|k|\ge s$,
 $$
 \frac 1{2\pi}\int_0^{2\pi}S^{[s]}_{N_j}(t)\mathrm{e}^{- \mathrm{i}kt}{dt}=
 \lambda_{|k|,r}(1-2^{-N_j})\frac{|k|!}{(|k|-s)!}\widehat f_k \mathop{\longrightarrow}\limits_{j\to\infty} \frac{|k|!}{(|k|-s)!}\widehat f_k.
 $$
Therefore, the equality  $g=f^{[s]}$ is true. Hence, for the function  $f$ and all $s=0,1,\ldots,r-n$, there exists the derivative $f^{[s]}$ and $f^{[s]}\in L_p$.

%%%Further
Now,
let us prove the estimate  (\ref{est Mp}). By virtue of (\ref{M_p_diff}), (\ref{OA}), for any $k\in\mathbb N$ and $\varrho\in (0,1)$, we have
\begin{equation}\label{deviation for AA}
M_p\left(\varrho, A_{k}-A_{k-1},r\right)\le \left\|A^{[r]}_{k}-A^{[r]}_{k-1}\right\|_p
=O\left(\frac{\omega(1-\varrho_{k-1})}{(1-\varrho_{k})^n}\right)
$$
$$
+O\left(\frac{\omega(1-\varrho_{k})}{(1-\varrho_{k-1})^n}\right)
=O\left(2^{kn}\omega(2^{-k+1})+2^{(k-1)n}\omega(2^{-k})\right)
$$
$$
=
O\left(2^{(k-1)n}\omega(2^{-(k-1)})\right),\ \  k\to +\infty.
\end{equation}
According to  (\ref{|f^[r]|}) and (\ref{f-Ap est}),  for any $r\in {\mathbb N}$, $\varrho\in (0,1)$ and
%%%$x\in[0, 2\pi]$
$x\in\mathbb T$, we obtain
\[%\begin{equation}\label{AA1}
M_p\left(\varrho, f-A_{\varrho,r}(f),r\right)\le 2r!\frac{\left\|f-A_{\varrho,r}(f)\right\|_{p}}{(1-\varrho)^{r}}
=O\left(\frac{\omega(1-\varrho)}{(1-\varrho)^{n}}\right), \quad \varrho\to 1-.
\]%\end{equation}
Therefore, for any positive integer $N$,
\begin{eqnarray}\nonumber
M_p\left(\varrho_{_{\scriptstyle N}}, f-A_{N}(f),r\right)&=& O\left(\frac{\omega(1-\varrho_{_{\scriptstyle N}})}{(1-\varrho_{_{\scriptstyle N}})^n}\right)\\ \label{deviation for f}
&=& O\left(2^{Nn}\omega(2^{-N})\right),\ N\to +\infty.
\end{eqnarray}

Consider the sum $\sum_{k=1}^N2^{(k-1)n}\omega(2^{-(k-1)})$, $N\in {\mathbb N}$.  Since the function  $\omega$ satisfies the condition $({\mathscr Z}_n)$ ,
%%%then
the function $\omega(t)/t^n$ almost decreases on $[0,1]$ (see, for example [\ref{Bari_Stechkin}]). Therefore,
\begin{equation}\label{summe}
 C_1\sum_{k=1}^N2^{(k-1)n}\omega(2^{-(k-1)}) \le  2^{(N-1)n}\omega(2^{-(N-1)})+\int_1^{N} 2^{(t-1)n}\omega(2^{-(t-1)})dt
 $$
 $$
  \le 2^{(N-1)n}\omega(2^{-(N-1)})+\frac 1{\ln 2}\int\limits_{2^{-N+1}}^1\! \! \! \omega(\tau)/\tau^{n+1}d\tau\le C_2 2^{(N-1)n}\omega(2^{-(N-1)}).
\end{equation}
Putting $\varrho=\varrho_{_{\scriptstyle N}}$ and taking into account relations (\ref{deviation for AA}), (\ref{deviation for f}), (\ref{summe}) and
$$
A_{0}(x)=S_{r-1}(f)(x)=\sum_{|k|\le r-1}\widehat f_k\mathrm{e}^{ \mathrm{i}kx},
$$
we get
\begin{equation}\label{M_p(varrho_N)}
M_p\left(\varrho_N,f,r\right)=M_p\left(\varrho_N,f-S_{r-1}(f),r\right)
$$
$$
=M_p\left(\varrho_N,f-A_{\varrho_{_{\scriptstyle N}}}+\sum\limits_{k=1}^N (A_{k}-A_{k-1}),r\right)
= O\left(\sum_{k=1}^{N}2^{(k-1)n}\omega(2^{-(k-1)})\right)
$$
$$
=O\left(2^{Nn}\omega(2^{-N})\right)
 =O\left(\frac{\omega(1-\varrho_{_{\scriptstyle N}})}{(1-\varrho_{_{\scriptstyle N}})^n}\right),\quad  N\to +\infty.
\end{equation}

If the  function  $\omega$ satisfies the condition $({\mathscr Z}_n)$ , then for all $t\in [0,1]$\ \  $\omega(2t)\le C \omega(t)$ (see, for example [\ref{Bari_Stechkin}]).  Furthermore, for all $\varrho\in[\varrho_{_{\scriptstyle N-1}},\varrho_{_{\scriptstyle
N}}]$, we have $1-\varrho_{_{\scriptstyle N}}\le
1-\varrho\le2(1-\varrho_{_{\scriptstyle N}})$. Hence,  relation (\ref{M_p(varrho_N)}) yields the estimate (\ref{est Mp}).

Now, applying the second inequality in Lemma 2 to the function   $f^{[r-n]}$, we get
\begin{eqnarray}\nonumber
K_{n}\left(1-\varrho, f^{[r-n]}\right)_p&\le& \|f^{[r-n]}-A_{\varrho, n}(f^{[r-n]})\|_p\\ \label{est Kn}
&+& \frac{4^{n}-1}3(1-\varrho)^nM_p(\sqrt{\varrho},f,r).
\end{eqnarray}

%%%%%%%%%%%%%%%%%%%%%%%%%%%%%%%%%%%%%%%%%%%%%%%%%%%%%%%%%%%%%
By virtue of (\ref{M_p_diff}) and (\ref{est Mp}), we see  that for  $\varrho\in [1/2,1)$,
\[
\int_\varrho^1 \bigg\|\frac{\partial^n f^{[r-n]}(\zeta,\cdot)}{\partial\zeta^n}\bigg\|_p(1-\zeta)^{n-1}d\zeta=\int_\varrho^1 \Big\|(f(\zeta,\cdot))^{[r]}(x)\Big\|_p\frac{(1-\zeta)^{n-1}}{\zeta^n} d\zeta
\]
\[
=\int_\varrho^1 M_p\left(\zeta, f,r\right)\frac{(1-\zeta)^{n-1}}{\zeta^n} d\zeta
\le 2^{n}C
\int_\varrho^1\frac{\omega(1-\zeta)}{1-\zeta}d\zeta
\]
\begin{equation}\label{111}
=O\left(\omega(1-\varrho)\right), \quad \varrho\to 1-.
\end{equation}
Therefore, we can apply  Lemma 4 to the function   $f^{[r-n]}$. Taking into account (\ref{M_p_diff}), we obtain
\[
f^{[r-n]}(x)-A_{\varrho, n}(f^{[r-n]})(x)=\frac{1}{(n-1)!}\int_\varrho^1 (f(\zeta,\cdot))^{[r]}(x)\frac{(1-\zeta)^{n-1}}{\zeta^n} d\zeta.
\]
Using   the integral Minkowski inequality and (\ref{111}), we %%%obtain
conclude
\begin{eqnarray}\nonumber
\|f^{[r-n]}-A_{\varrho, n}(f^{[r-n]})\|_p&\le& \frac{1}{(n-1)!}\int_\varrho^1 M_p\left(\zeta, f,r\right)\frac{(1-\zeta)^{n-1}}{\zeta^n} d\zeta\\ \label{112}
&=& O\left(\omega(1-\varrho)\right), \quad \varrho\to 1-.
\end{eqnarray}
Combining relations (\ref{est Kn}), (\ref{est Mp}) and (\ref{112}),  we get (\ref{K-funct est}).

%%%%%%%%%%%%%%%%%%%%%%%%%%%%%%%%%%%%%%%%%%%%%%%%%%%%%%%%%%%%%

\vskip 0.5cm

{\bf Acknowledgments.} This work was supported in part by the FP7-People-2011-IRSES project number  295164 (EUMLS:
%%%EUUkrainian
 EU-Ukrainian
 Mathematicians for Life Sciences).

\vskip 3mm

%%%%%%%%%%%%%%%%%%%%%%%%%%%%%%%%%%%%%%%%%%%%%%%%%%%%%%%%%%%%%%%%%%%%%%%%%%%%%%%%%%%%%%%%%%%%%%%%%%%%%%%%%%%%%%%%%%%%%%%%%%%%%%%%%%%%%%%%%%%%%%%%%%%%%%%%%

\small

\begin{enumerate}

%%%%%%%%%%%%%%%%%%%%%%%%%%%%%%%%%%%%%%%%%%%%%%%%%%%%%%%%%%%%%%%%%%%%%%%%%%%%%%%%%%%%%%%%%%%%%%%%%%%%%%%%%%%

  \item\label{Bari_Stechkin} %15
 {\it Bari  N.\,K., Stechkin S.\,B.\/}  Best approximations and differential properties of two conjugate functions // Tr. Mosk. Mat. Obshch. --- 1956.
 --- 5. --- P.~483--522 .

\item \label{Butzer_Nessel}
{\it Butzer P., Nessel R.\/} Fourier Analysis and Approximation.
One --- Dimentional Theory. --- Basel--New York, 1971. --- 554~p.

\item \label{Butzer_Sunouchi}
 {\it Butzer P.\,L., Sunouchi G.\/}  Approximation theorems for the solution of Fourier's problem and Dirichlet's problem  // Math. Ann. --- 1964.  --- 155. --- P.~316--330.

\item \label{Butzer_Tillmann}
 {\it Butzer P.\,L., Tillmann H.\,G.\/}  Approximation theorems for semi-groups of bounded linear transformations // Math. Ann. --- 1960.  --- 140. --- P.~256--262.

\item \label{Butzer}
 {\it Butzer P.\,L.\/}   Beziehungen zwischen den Riemannschen, Taylorschen und
 %%%gewohnlichen
 gew\"ohnlichen
 Ableitungen reellwertiger Funktionen // Math. Ann. --- 1961.  --- 144. --- P.~275--298.

 \item \label{Chandra}
 {\it  Chandra P., Mohapatra R.\,N.\/}  Approximation of functions by $(J,q_n)$ means of Fourier series // Approx. Theory Appl. --- 1988. --- 4, №~2. ---  P.~49--54.

\item \label{Chui}
 {\it Chui C.\,K.,  Holland A.\,S.\,B.\/}  On the order of approximation by Euler and Taylor means // J. Approx. Theory. --- 1983. --- 39, №~1. ---  P.~24--38.

 \item\label{DeVore Lorentz}
{\it DeVore R.\,A.,  Lorentz G.\,G.\/}  Constructive approximation. --- Berlin: Springer--Verlag, 1993. --- 449~p.

\item \label{Holland}
 {\it Holland A.\,S.\,B., Sahney B.\,N., Mohapatra R.\,N.\/}  $L_p$ approximation of functions by Euler means // Rend. Mat. --- 1983. --- 3 (7), №~2. ---  P.~341--355.

\item \label{Leis}
 {\it Leis R.\/}  Approximationss{\"{a}}tze f{\"{u}}r stetige Operatoren // Arch. Math. --- 1963. --- 14. ---  P.~120--129.

\item \label{Mohapatra}
 {\it Mohapatra R.\,N.,  Holland A.\,S.\,B.,  Sahney B.\,N.\/}  Functions of class ${\rm Lip}(\alpha,p)$ and their Taylor mean // J. Approx. Theory. --- 1985. --- 45, №~4. ---  P.~363-–374.

\item \label{Prestin_Savchuk_Shidlich}
 {\it Prestin J., Savchuk V.\,V., Shidlich A.\,L.\/} Approximation  of  $2\pi$-periodic functions by Taylor--Abel--Poisson operators in the integral metric // to appear in Dopov. NANU. --- 2017. --- №~1.

\item \label{Rudin}
 {\it  Rudin W.\/}  Function theory in polydiscs. --- Moscow: Mir,  1974. --- 160~p.

\item \label{Savchuk}
 {\it Savchuk V.\,V.\/}  Approximation of holomorphic functions by Taylor-Abel-Poisson means //
 Ukr. Mat. Zh. --- 2007. --- 59, №~9. --- P.~1253--1260.

\item \label{Savchuk_Shidlich}
 {\it Savchuk V.\,V.,  Shidlich A.\,L.\/}  Approximation of functions of several variables by linear methods in the space $S^p$ // Acta Sci. Math. --- 2014. --- 80, №~3--4. --- P.~477--489.

\item \label{Trigub_Bellinsky}
 {\it Trigub R.\,M.,  Bellinsky E.\,S.\/}  Fourier analysis and approximation of functions. ---  Dordrecht: Kluwer Academic Publishers, 2004. --- 585~p.

\item \label{Savchuk_Zastavnyi}
 {\it Zastavnyi V.\,P., Savchuk V.\,V.\/} Approximation of classes of convolutions by linear operators of a special form // Mat. Zametki. --- 2011. --- 90, №~3. --- P.~351--361.

\end{enumerate}

\enddocument